\documentclass[12 pt]{amsart}
\usepackage{amsmath,mathrsfs,amsfonts,amssymb,xypic,fullpage,setspace, stmaryrd, enumerate, geometry,hyperref }

\usepackage[all]{xy}

\DeclareFontFamily{U}{wncy}{}
\DeclareFontShape{U}{wncy}{m}{n}{<->wncyr10}{}
\DeclareSymbolFont{mcy}{U}{wncy}{m}{n}
\DeclareMathSymbol{\Sh}{\mathord}{mcy}{"58} 

\theoremstyle{plain}
\newtheorem{thm}{Theorem}
\newtheorem{cor}{Corollary}
\newtheorem{lem}{Lemma}
\newtheorem{prop}{Proposition}

\newenvironment{pf}[1][Proof]{\begin{trivlist}
\item[\hskip \labelsep {\bfseries #1}]}{\end{trivlist}}

\newenvironment{example}[1][Example]{\begin{trivlist}
\item[\hskip \labelsep {\bfseries #1}]}{\end{trivlist}}
\newenvironment{rmk}[1][Remark]{\begin{trivlist}
\item[\hskip \labelsep {\bfseries #1}]}{\end{trivlist}}

\newcommand{\rhobar}{\bar{\rho}}
\newcommand{\epad}{\epb \otimes \mathrm{ad} \ \bar{\rho}_{f,\lambda}}
\newcommand{\epadzero}{\epb \otimes \mathrm{ad}^0 \ \bar{\rho}_{f,\lambda}}
\newcommand{\epadzerop}{\epb \otimes \mathrm{ad}^0 (\bar{\rho}_{f,\lambda})}
\newcommand{\Q}{\mathbf{Q}}

\newcommand{\Z}{\mathbf{Z}}

\newcommand{\ep}{\epsilon_\ell}
\newcommand{\epb}{\bar{\epsilon}_\ell}

\newcommand{\barFl}{\bar{\mathbf{F}}_\ell}

\newcommand{\Qp}{\mathbf{Q}_p}

\title{Obstruction Criteria for Modular Deformation Problems}
\author{Jeffrey Hatley}
\date{\today}

\begin{document}
\maketitle


\section{Introduction}
Let $f=\sum a_n q^n$ be a newform of weight $k \geq 2$ and level $\Gamma_1(N)$. Write $K=\Q(a_n)$ for the number field generated by its Fourier coefficients. Let $\lambda$ be a prime of $K$, and let $\ell$ be the characteristic of its residue field $k_\lambda$. For any finite set $S$ of places which contains the primes dividing $N \infty$, let $\Q_{S \cup \{\ell\}}$ be the maximal extension of $\Q$ unramified outside $S \cup \{\ell\}$, and let $G_{\Q,S \cup \{\ell\}}$ be its Galois group over $\Q$. By work of Deligne, there is an associated semisimple residual Galois representation
\[
\bar{\rho}_{f,\lambda} : G_{\Q,S \cup \{\ell\}} \to \mathrm{GL}_2 (k_\lambda),
\]
and this representation is absolutely irreducible for almost all primes $\lambda$. 

Given such a representation $\rhobar$, it is interesting to study its lifts to other coefficient rings. If $A$ is a local ring with residue field $k_\lambda$, we say $\rho$ is a \textit{lift} of $\rhobar$ if the following diagram commutes:
\begin{displaymath}
    \xymatrix{ G_\Q \ar[r]^{\rho}  \ar[dr]_{\rhobar_{f,\lambda}} & \mathrm{GL}_2 (A) \ar[d] \\
                &  \mathrm{GL}_2 (k_\lambda) }
\end{displaymath}
\noindent The vertical arrow is induced by the reduction map $A \to k_\lambda$; we consider two lifts equivalent if they are conjugate to one another by a matrix in the kernel of this induced map. An equivalence class of lifts is called a \textit{deformation} of $\rhobar$.

The study of the deformation theory of such Galois representations, which began with Mazur's seminal paper \cite{Maz}, has been the subject of much important research in number theory; in particular, it featured prominently in the proof of the Taniyama-Shimura conjecture, and more recently, in the proof of Serre's Conjecture. See Section \ref{Primer} for a brief introduction to deformation theory and the terms used below.

In a pair of papers \cite{WestonUnob, WestonEx}, Weston proved that for any newform $f$ of weight $k \geq 2$, the deformation problem for $\rhobar_{f,\lambda}$ is unobstructed for infinitely many primes $\lambda$, and when the level of $f$ is squarefree, he gave an explicit description of the obstructed primes. In fact, when $k \geq 3$, there are only finitely many obstructed primes, while for $k =2$ the obstructed primes are a set of density zero. The first main result of this paper is the removal of the squarefree hypothesis from Weston's result; only a minor modification of the bound given in \cite{WestonEx} is necessary. See Theorem \ref{MainThm1Body} in Section \ref{Squarefree} for the full statement.

While Theorem \ref{MainThm1Body} gives sufficient conditions for a deformation problem to be unobstructed, the second main result of this paper focuses on a necessary condition. For any modular Galois representation $\bar{\rho}: G_\Q \to \mathrm{GL}_2 (\barFl)$, there is an optimal (least) level $N$ coprime to $\ell$ such that $\bar{\rho}$ arises from a newform of level $N$. Call a deformation problem \textit{minimal} if the set $S$ of primes (as in the first paragraph) contains only those places dividing $N\infty$. We show in Theorem \ref{MainThm2} that minimal deformation problems are only unobstructed when they arise from modular forms of optimal level. This is analogous to a similar phenomenon which occurs in Hida Hecke algebras.

\textbf{Notation.}
We fix an algebraic closure $\bar{\Q}$ of $\Q$, and for each rational prime $\ell$, we fix an embedding $\bar{\Q} \hookrightarrow \bar{\Q}_\ell$. Let $G_\Q=$Gal$(\bar{\Q}/\Q)$ and let $G_\ell$=Gal$(\bar{\Q}_\ell/\Q_\ell)$. Whenever $S$ is a finite set of primes, $G_{\Q,S}$ denotes the Galois group (over $\Q$) of the maximal extension of $\Q$ which is unramified outside of $S$.

We write $\ep$ for the $\ell$-adic cyclotomic character. For any character $\psi$ we denote its reduction mod $\lambda$ by $\bar{\psi}$, where $\lambda$ is made clear in context.

If $\rho : G \to V$ is a representation, the adjoint representation $\mathrm{ad}\ \rho : G \to \mathrm{End}(V)$ is defined by letting $g \in G$ act on End$(V)$ via conjugation by $\rho(g)$; we write ad$^0 \rho$ for the trace-zero component of the adjoint.

\section{Deformation Theory}\label{Primer}

Consider an odd, continuous Galois representation $\bar{\rho}:G_{\Q,S} \to \mathrm{GL}_2 (\mathbf{F})$, where $\mathbf{F}$ is some finite field and $S$ is a finite set of primes containing the characteristic of $\mathbf{F}$ and the infinite place. Let $\mathcal{C}$ be the category whose objects are local rings which are inverse limits of artinian local rings with residue field $\mathbf{F}$, and whose morphisms $A \to B$ are continuous local homomorphisms inducing the identity map on residue fields. As explained in the introduction, if $A \in \mathcal{C}$, then we say $\rho: G_{\Q,S} \to \mathrm{GL}_2 (A)$ is a lift of $\rhobar$ if the composition
\[
G_{\Q,S} \xrightarrow{\rho} \mathrm{GL}_2 (A) \to \mathrm{GL}_2 (\mathbf{F})
\]
is equal to $\rhobar$. Two lifts $\rho_1,\rho_2$ of $\bar{\rho}$ to $A$ are considered equivalent if they are conjugate to one another by a matrix in the kernel of the map $\mathrm{GL}_2 (A) \to \mathrm{GL}_2 (\mathbf{F})$, and a \textit{deformation} of $\rhobar$ to $A$ is an equivalence class of lifts of $\rhobar$ to $A$. There is an associated \textit{deformation functor}
\[
D_{\bar{\rho}}^S : \mathcal{C} \to \mathrm{Sets}
\]
which sends a ring $A$ to the set of deformations of $\bar{\rho}$ to $A$. When $\bar{\rho}$ is absolutely irreducible, this functor is representable by a ring $R_{\bar{\rho}} \in \mathcal{C}$ \cite{Maz,MDic}. 

For $i=1,2$, let $d_i$ be the $\mathbf{F}$-dimension of the Galois cohomology group $H^i(G_{\Q,S}, \mathrm{ad} \ \bar{\rho})$. Mazur showed that that $d_1 - d_2 \geq 3$ and
\[
R_{\bar{\rho}} \simeq W(\mathbf{F}) \llbracket T_1, \ldots, T_{d_1} \rrbracket / (r_1, \ldots, r_{d_2}),
\]
where $W(\mathbf{F})$ is the ring of Witt vectors of $\mathbf{F}$. When $d_2=0$, it can be shown that $d_1=3$, so $R_{\bar{\rho}}$ is simply a power series ring in three variables. In this case, the deformation problem for $\bar{\rho}$ is said to be \textit{unobstructed}.

Let $\ell$ be the characteristic of $\mathbf{F}$. As in Lemma 2.5 of \cite{WestonEx}, an application of the Poitou-Tate exact sequence allows one to show that
\begin{equation}\label{h2dim}
\mathrm{dim}_\mathbf{F} H^2(G_{\Q,S},\mathrm{ad} \ \bar{\rho}) \leq \mathrm{dim}_\mathbf{F} \Sh^1(G_{\Q,S}, \epb \otimes \mathrm{ad}^0   \bar{\rho}) + \sum_{p \in S} \mathrm{dim}_\mathbf{F} H^0(G_p, \epb \otimes \mathrm{ad} \ \bar{\rho})
\end{equation}
with equality if $\ell \neq 3$ . Here $\Sh^1(G_{\Q,S}, \epb \otimes \mathrm{ad}^0   \bar{\rho})$ is a sort of Selmer group; when $\bar{\rho}=\bar{\rho}_{f,\lambda}$ for some newform $f$, this term can be controlled by the set Cong$(f)$ of congruence primes for $f$, as described in (\cite{WestonEx}, Section 4). Our focus will instead be on the local invariants $H^0(G_{\Q,S}, \epb \otimes \mathrm{ad} \ \bar{\rho})$ for $p \in S$, which we refer to as \textit{obstructions at $p$}. 

\section{Removing the squarefree hypothesis}\label{Squarefree}

We fix some notation to be used throughout Section \ref{Squarefree}. Let $f=\sum a_n q^n$ be a newform of level $N$ and weight $k \geq 2$. Let $\omega$ be its nebentypus character, and let $M$ be the conductor of $\omega$. Let $K$ be its associated number field, and fix a prime $\lambda$ in $K$ with residue field $k_\lambda$ of characteristic $\ell$ such that $(N,\ell)=1$ and $\bar{\rho}_{f,\lambda}$ is absolutely irreducible. Let $S$ be a finite set of places containing the primes dividing $N\infty$. We wish to study the conditions under which the deformation problem for
\[
\bar{\rho}_{f,\lambda} : G_{\Q,S \cup \{\ell \}} \to \mathrm{GL}_2 (k_\lambda)
\] 
is unobstructed, and as described in Section \ref{Primer}, as long as $\lambda \notin \mathrm{Cong}(f)$, then this amounts to determining when $H^0(G_p, \epb \otimes \mathrm{ad} \ \bar{\rho}) \neq 0$ for $p \in S$. 

Let $\pi$ be the automorphic representation associated to $f$, and write $\pi = \otimes' \pi_p$ for its decomposition into admissible complex representations $\pi_p$ of $\mathrm{GL}_2 (\Qp)$.  By the Local Langlands correspondence, the classification of each $\pi_p$ allows us to study $\bar{\rho}_{f,\lambda}|_{G_p}$ in an explicit fashion. In \cite{WestonEx}, the assumption that $N$ be squarefree aided in the determination of $\pi_{p}$ for each $p \in S$; in particular, in this case it is easy to determine when $\pi_{p}$ is an unramified principal series, a principal series with one ramified character and one unramified character, or a special (twist of Steinberg) representation, and these are the only possibilities. When $p^2 \mid N$, it is not so easy to determine the structure of $\pi_{p}$.  However, determining the exact structure of $\pi_p$ turns out to be unnecessary. 

\subsection{Twists and $p$-primitive newforms} Recall that for any primitive Dirichlet character $\chi$ of conductor $M$, we may twist the newform $f$ to obtain a newform $f \otimes \chi=\sum b_n q^n$, where $b_n = \chi(n) a_n$ for almost all $n$. The level of $f \otimes \chi$ is at most $NM^2$, but it may be smaller. For any newform $f$ and any prime $p$, one says that $f$ is \textit{$p$-primitive} if the $p$-part of its level is minimal among all its twists by Dirichlet characters. We have the following lemma.

\begin{lem}\label{adtwist} Let $f$ be a newform and let $f_p$ be a $p$-primitive twist. Then
\[
H^0(G_p,\epad) = H^0(G_p,\epb \otimes \mathrm{ad} \ \rhobar_{f_p, \lambda}). 
\]
In particular, $f$ has local obstructions at $p$ if and only if $f_p$ has local obstructions at $p$.
\end{lem}
\begin{pf}
For some Dirichlet character $\chi$ we have $f_p = f \otimes \chi$. It follows that $\rhobar_{f_p,\lambda} \simeq \chi \otimes \rhobar_{f,\lambda}$, and a straightforward matrix calculation then shows that ad$(\rhobar_{f_p,\lambda}) \simeq $ ad($\rhobar_{f,\lambda}$). The lemma follows.
\qed \end{pf}

By Lemma \ref{adtwist}, when studying local obstructions at $p$ for a newform $f$, we may assume that $f$ is $p$-primitive.
The utility of considering $p$-primitive newforms is given by the following result, which comes from  (\cite{LW}, Proposition 2.8):

\begin{prop}\label{LWprop}
Let $\pi_p$ be the local component of a $p$-primitive newform $f \in S_k (\Gamma_1(Np^r))$ with $p \nmid N$ and $r \geq 1$. Then one of the following holds.
\begin{enumerate}[(1)]
\item $\pi_p \simeq \pi(\chi_1,\chi_2)$ is principal series, where $\chi_1$ is unramified and $\chi_2$ is ramified;
\item $\pi_p \simeq \mathrm{St} \otimes \chi$, is special (twist of Steinberg) with $\chi$ unramified;
\item $\pi_p$ is supercuspidal.
\end{enumerate}
\end{prop}

\begin{pf}
See  (\cite{LW}, Proposition 2.8) for the proof.
\qed \end{pf}

\begin{rmk}If the level of a newform $f$ is divisible by $p^2$, it may be difficult to explicitly determine its $p$-minimal twist. Loeffler and Weinstein have made this computationally feasible in many cases; see \cite{LW}. We will avoid this extra difficulty and simply determine where obstructions might occur in all three cases of the above proposition.\end{rmk}

\subsection{Supercuspidal Obstruction Conditions}

The arguments used by Weston in \cite{WestonEx} are robust enough to carry over into the non-squarefree setting when we are cases (1) and (2) of Proposition \ref{LWprop}. We instead focus on case (3), where $\pi_p$ is supercuspidal. We will frequently make use of the fact that
\[
\mathrm{dim}_k H^0(G_p, \epb \otimes (\mathrm{ad}^0 \rhobar_{f,\lambda})^{\mathrm{ss}} ) \leq \mathrm{dim}_k H^0(G_p,\epadzero).
\]

When $p > 2$, a supercuspidal $\pi_p$ is always induced from a quadratic extension of $\Q_p$, and these will be the focus of Proposition \ref{SCProp} below. When $p=2$, there are additional supercuspidal representations, called extraordinary representations, and we consider these first. The case where $\pi_{p}$ is extraordinary was actually already dealt with in (\cite{WestonUnob}, Proposition 3.2) and are not a problem if $\ell \geq 5$. We reproduce the proof here.

\begin{prop}
Suppose $\pi_p$ is extraordinary, so $p=2$, Then $H^0(G_p, \epad) = 0$ if $\lambda$ has residue characteristic at least $5$.
\end{prop}
\begin{pf}
Let $\rho: G_2 \to \mathrm{GL}_2 (\bar{\Q}_\ell)$ be the representation of $G_2$ which is in Langlands correspondence with $\pi_2$. In this case, the projective image of inertia,  proj $\rho(I_2)$, in $\mathrm{PGL}_2 (\bar{\Q}_\ell)$ is isomorphic to either $A_4$ or $S_4$, and the composition
\[
\mathrm{proj} \  \rho (I_2) \hookrightarrow \mathrm{PGL}_2 (\bar{\Q}_\ell) \xrightarrow{\mathrm{ad}^0} \mathrm{GL}_3 (\bar{\Q}_{\ell})
\]
is an irreducible representation of $\mathrm{proj} \ \rho$. Since $\mathrm {proj} \ \rho(I_2)$ has order $12$ or $24$, it follows that $\mathrm{ad}^0 \bar{\rho}_{f,\lambda}$ is an irreducible $\bar{\mathbf{F}}_\ell$-representation of $I_2$ since char$(\lambda) \geq 5$, thus $H^0(I_2,\epadzero)=0$ and the proposition follows.
\qed \end{pf}

We will henceforth assume $\ell \geq 5$, so by this proposition there are no obstructions in the extraordinary case.

Now we deal with the final remaining possibility for $\pi_{p}$, which is the supercuspidal, non-extraordinary case. Recall that for any character $\psi$, we write $\bar{\psi}$ for its reduction mod $\lambda$. 

\begin{prop}\label{SCProp}
Suppose $f$ is a newform of weight $k \geq 2$ such that $\pi_{p}$ is supercuspidal but not extraordinary. Suppose also that $\ell > 5$. If $p^4 \not\equiv 1$ (mod $\ell$), then $H^0(G_p, \epad) = 0$.
\end{prop}

\begin{pf}
The Langlands correspondence (cf. \cite{WestonUnob}, Proposition 3.2 or \cite{LW}, Remark 3.11) implies that there is a quadratic extension $E/\Qp$ such that in characteristic zero we have
\[
\rho_{f,\lambda}|_{G_p} \simeq \mathrm{Ind}_{G_E}^{G_p} \chi
\]
where $G_E=$Gal$(\bar{E}/E)$ is the absolute Galois group of $E$ and $\chi : G_E \to \bar{\Q}_\ell$ is a continuous character. Let $\chi_E : \mathrm{Gal}(E/\Qp) \to \{ \pm 1 \}$ be the nontrivial character for $E/ \Qp$. Let $\chi^c$ be the Galois conjugate character of $\chi$, and let $\psi=\chi \cdot (\chi^c)^{-1}$. We have
\[
\epadzerop^{\mathrm{ss}} \simeq \epb \chi_E \ \oplus \ \left( \epb \otimes \mathrm{Ind}_{G_E}^{G_p} \ \bar{\psi} \right)
\]
Since $\ell >3$, the first summand has no $G_p$-invariants, so we may focus on the second summand. By Mackey's criterion, the induced representation $\mathrm{Ind}_{G_E}^{G_p} \ \bar{\psi}$ is irreducible if and only if $\bar{\psi} \neq \bar{\psi^c}$. If it is irreducible, then so is its twist and we are done. 

So suppose that $\bar{\psi} = \bar{\psi^c}$. We first note that, since $\bar{\psi}=\bar{\chi} (\bar{\chi}^c)^{-1}$, we have $\bar{\psi}^c=\bar{\chi}^c \bar{\chi}^{-1}$, hence 
\[
\bar{\psi}^2=\bar{\psi}\bar{\psi}^c=[\bar{\chi} (\bar{\chi}^c)^{-1}]\cdot[\bar{\chi}^c (\bar{\chi})^{-1}]=1.
\]
Thus, $\bar{\psi}$ is a quadratic character on $G_E$. 

Restricting the induced representation to $G_E$ we have
\[
(\mathrm{Ind}_{G_E}^{G_p} \bar{\psi} )|_{G_E} \simeq \bar{\psi} \oplus \bar{\psi}^c = \bar{\psi} \oplus \bar{\psi}
\] 
where the first equality is a generality about induced representations and the second comes from our assumption that $\bar{\psi}=\bar{\psi}^c$. So we already have $H^0(G_E,\epb \otimes \mathrm{Ind}_{G_E}^{G_p}\bar{\psi})=0$ unless $\bar{\psi}=\epb|^{-1}_{G_E}$, in which case $\epb|_{G_E}$ is quadratic. Since $G_E$ has index $2$ in $G_p$, this would imply that on $G_p$ the cyclotomic character $\epb$ has order at most $4$. Evaluating at $\mathrm{Frob}_p$, this implies that $p^4 \equiv 1$ (mod $\ell$). So if $p^4 \not\equiv 1$ (mod $\ell$), then the representation has no $G_E$-invariants and hence it has no $G_p$ invariants, completing the proof.
\qed \end{pf}

We are now ready to prove the first main theorem, which removes the squarefree hypothesis from (\cite{WestonEx}, Theorem 4.3). For a newform $f$ of level $N$, Cong$(f)$ is the set of congruence primes for $f$, i.e. the primes $\lambda$ such that there exists a newform $g$ (which is not a Galois conjugate of $f$) of level dividing $N$ with $f \equiv g$ (mod $\lambda$).

\begin{thm}\label{MainThm1Body}
Assume that $\bar{\rho}_{f,\lambda}$ is absolutely irreducible and $\ell > 3$. If $H^2(G_{\Q,S \cup \{\ell\}}, \epad) \neq 0$ then one of the following holds:
\begin{enumerate}
\item $\ell \leq k$;
\item $\ell \mid N$;
\item $\ell \mid \phi(N_S)$, where $N_S$ is the product of the primes in $S$;
\item $\ell \mid (p+1)$ for some $p \mid N$;
\item $a_p^2 \equiv (p+1)^2 p^{k-2} \omega(p)$ (mod $\lambda$) for some $p \in S$, $p \nmid N$, $p \neq \ell$;
\item $\ell = k+1$ and $f$ is ordinary at $\lambda$;
\item $k=2$ and $a_\ell^2 \equiv \omega(\ell)$ (mod $\lambda$);
\item $N=1$ and $\ell \mid (2k-3)(2k-1)$;
\item $\lambda \in $ Cong$(f)$;
\item $p^4 \equiv 1$ (mod $\ell$) for some $p$ such that $p^2 \mid N$.
\end{enumerate}
\end{thm}

\begin{rmk}We note that conditions (1)--(9) are essentially the same conditions from (\cite{WestonEx}, Theorem 4.3); these conditions deal with the non-supercuspidal primes in $S$, while condition (10) deals with the (potentially) supercuspidal primes.\end{rmk}

\begin{pf}
By equation (\ref{h2dim}), if $H^2(G_{\Q,S \cup \{\ell\}}, \epad) \neq 0$ then either $\mathrm{dim}_k \Sh^1(G_{\Q,S}, \epb \otimes \mathrm{ad}^0   \bar{\rho}) \neq 0$ or $H^0(G_p,\epad)\neq0$ for some $p \in S$. By (\cite{WestonEx}, Lemma 17), the former is only possible if $\lambda \in$ Cong$(f)$. This is accounted for in condition (9).

Now let $p \in S$. While determining whether $H^0(G_p,\epad)=0$, Lemma \ref{adtwist} allows us to replace $f$ by its $p$-minimal twist. In this case, by Lemma \ref{LWprop} there are only three possibilities for the local representation $\pi_p$.

If $\pi_{p}$ is principal series or special as in cases (1) and (2) of Lemma \ref{LWprop}, then the local Galois representation has exactly the same form as the cases handled by Weston (\cite{WestonEx}, Theorem 4.3). This accounts for conditions (1)--(9). The only difference occurs in condition (4). In Weston's original condition, it is only necessary to avoid $\ell \mid (p+1)$ for primes $p$ dividing $N/M$, where $M$ is the conductor of the nebentypus character of $f$. Since we have replaced $f$ by its $p$-minimal twist $f_p$, and we don't know the conductor of the character of $f_p$, we replace Weston's original condition with our coarser condition.

If $\pi_{f,p}$ is supercuspidal, then Proposition \ref{SCProp} yields condition (10). This covers all the possibilities for $\pi_{f,p}$, thus completing the proof. 
\qed \end{pf}

\section{Minimal Deformation Problems and Optimal Levels}\label{Minimal}

Given a modular form $f$, a prime $\lambda$ of $\bar{K}$, and a finite set of places $S$, let us write $\mathbf{D}(f,S)$ for the corresponding deformation problem. (We suppress $\lambda$ from the notation, as it will always be clear from context.) If $S$ contains only the primes dividing the level of $f$ and the infinite place, then we may simply write $\mathbf{D}(f)$, and we call this the \textit{minimal deformation problem} for $f$. 

For any odd, continuous, absolutely irreducible representation $\bar{\rho} : G_\Q \to \mathrm{GL}_2(k_\lambda)$, with $k_\lambda$ a finite field of characteristic $\ell$, let $\mathcal{H}(\bar{\rho})$ be the set of newforms of level prime to $\ell$ giving rise to this representation. Among all such newforms, there is a least level appearing, which we call the \textit{optimal} level for $\mathcal{H}(\bar{\rho})$. In fact, this optimal level is the prime-to-$\ell$ Artin conductor of $\bar{\rho}$ (see \cite{DT}).

Let $f$ and $g$ be newforms in $\mathcal{H}(\bar{\rho})$ with associated minimal sets of primes $S$ and $S'$, respectively. We have an isomorphism of residual Galois representations $
\rho_{f,\lambda} \simeq \rho_{g,\lambda}$, and if $S \subset S'$ then we have an equality of deformation problems $\mathbf{D}(f,S')=\mathbf{D}(g)$. Furthermore, since $S \subset S'$, if $\mathbf{D}(f)$ is obstructed then so is $\mathbf{D}(g)$. In fact, we prove the following theorem:

\begin{thm}\label{MainThm2}
If $\mathbf{D}(f)$ is unobstructed, then $f$ is of optimal level for $\mathcal{H}(\bar{\rho})$.
\end{thm}

In Section 4.2 we present the proof of this theorem; our strategy is to prove the contrapositive. By Proposition \ref{Mform} below, we know the factorization of any nonoptimal level. If $g$ is a newform of nonoptimal level, we compare it to an optimal level newform $f$. Since, as discussed above, $\mathbf{D}(g)$ inherits any obstructions that $\mathbf{D}(f)$ might have, we may assume that $\mathbf{D}(f)$ is unobstructed, and we show that even in this case, $\mathbf{D}(g)$ is necessarily obstructed.

This theorem is motivated by the following heuristic: If $\bar{\rho}$ is $\ell$-ordinary and $\mathcal{H}(\bar{\rho})$ is a Hida family, then its components of non-optimal level have associated (full) Hecke algebras of higher $\Lambda = \Z_\ell \llbracket T \rrbracket$-rank than the optimal-level component (cf. \cite{EPW}, Section 2.4). Thus, if a general enough $\mathcal{R}=\mathbf{T}$ theorem is known (or believed), then this forces the deformation ring to grow as well. Our theorem shows that this sort of behavior is not a special property of Hida families, and that it actually occurs independent of any geometric structure.

\begin{rmk}
It is worth pointing out two things about this theorem. The first is that it does \textit{not} follow immediately from Theorem \ref{MainThm1Body}, because condition (9) of that theorem is not a sharp obstruction criterion, i.e. it does not \textit{guarantee} the existence of obstructions.  The other noteworthy aspect is that in \cite{WestonEx}, congruence primes are shown to (potentially) give rise to \textit{global} obstruction classes, whereas our proof uses the existence of a newform congruence to produce \textit{local} obstruction classes.
\end{rmk}

\subsection{Preliminaries}

In this section we record the results which we will use to prove the theorem. Let us first set some notation to be used throughout Section \ref{Minimal}.

Let $f = \sum a_n q^n$ be a newform of weight $k \geq 2$, level $N$ (coprime to $\ell$), and nebentypus $\omega$, and let $M$ be the conductor of $\omega$. Let $S$ be a finite set of places containing the primes which divide $N\infty$. Let $K=\Q(a_n)$, and fix a prime $\lambda$ of $\bar{K}$ which lies over $\ell$. We have $f \in \mathcal{H(\rhobar)}$, where $\rhobar_{f,\lambda} \simeq \rhobar$. 

Suppose $f$ is of optimal level for $\mathcal{H}(\rhobar)$. If $g \in \mathcal{H}(\rhobar)$ is of nonoptimal level, we will want to know what form its level can have. The following is a result of Carayol (see the introduction of \cite{DT}).

\begin{prop}\label{Mform}
Suppose $\rho : G_\Q \to \mathrm{GL}_2(\barFl)$ is modular of weight $k \geq 2$ and level $N'$ coprime to $\ell$. Then 
\[
N'=N \cdot \prod \ p^{\alpha(p)}
\]
where $N$ is the conductor of $\bar{\rho}$, and for each $p$ with $\alpha(p)>0$, one of the following holds:
\begin{enumerate}
\item $p \nmid N \ell$, $p(\mathrm{tr} \rho(\mathrm{Frob}_p)^2) = (1+p)^2 \mathrm{det} \rho(\mathrm{Frob}_p)$ in $\barFl$ and $\alpha(p)=1$;
\item $p \equiv -1$ mod $\ell$ and one of the following holds:
\begin{enumerate}
\item $p \nmid N$, tr $(\rho(\mathrm{Frob}_p)=0$ in $\barFl$ and $\alpha(p)=2$, or
\item $p \mid \mid N$, $\mathrm{det}$ $\rho$ is unramified at $p$, and $\alpha(p)=1$;
\end{enumerate}
\item $p \equiv 1$ mod $\ell$ and one of the following holds:
\begin{enumerate}
\item $p \nmid N$ and $\alpha(p)=2$, or
\item $p^2 \nmid N$, or the power of $p$ dividing $N$ is the same as the power dividing the conductor of det $\rho$, and $\alpha(p)=1$.
\end{enumerate}
\end{enumerate}
\end{prop}

Our goal, then, is to show that each of the possible supplementary primes appearing in Proposition \ref{Mform} gives rise to an obstruction. We collect some lemmas in this direction.

The first two lemmas come from \cite{WestonEx}.

\begin{lem}\label{LemmaPhi}
If $\ell \mid (p-1)$ for some $p \in S$, then $\mathbf{D}(f,S)$ is obstructed. 
\end{lem}

\begin{pf}
This is proved in the discussion at the beginning of Section 3 of \cite{WestonEx}. Since $\ell \mid p-1$, we have $H^0(G_p, \bar{\epsilon}_\ell) \neq 0$. Since $\bar{\epsilon}_\ell \otimes \mathrm{ad} \ \bar{\rho}_{g,\lambda} \simeq \bar{\epsilon}_\ell \oplus (\bar{\epsilon}_\ell \otimes \mathrm{ad}^0\bar{\rho}_{g,\lambda})$ this shows that $H^0(G_p,\bar{\epsilon}_\ell \otimes \mathrm{ad} \ \bar{\rho}_{g,\lambda}) \neq 0$ and so $\mathbf{D}(f,S)$ is obstructed.
\qed \end{pf}

The previous lemma gives us a tool we can use when $p \equiv 1$ (mod $\ell$). The next lemma deals with the case when $p \not\equiv 1$ (mod $\ell$) and $p \neq \ell$.

\begin{lem}\label{LemmaAP}
Assume $p \nmid N \ell$ and $p \not\equiv 1$ (mod $\ell$). Then $H^0(G_p,\bar{\epsilon}_\ell \otimes \mathrm{ad} \ \bar{\rho}_{f,\lambda}) \neq 0$ if and only if $a_p^2 \equiv (p+1)^2 p^{k-2} \omega(p)$ (mod $\lambda$).
\end{lem}

\begin{pf}
This is (\cite{WestonEx}, Lemma 3.1).
\qed \end{pf}

In this final lemma, we prove a partial converse to (\cite{WestonUnob}, Lemma 3.3).

\begin{lem}\label{LemmaSpecial}
If $p \mid \mid N$, $p \nmid M$, and $p^2 \equiv 1$ (mod $\ell$), then $\mathbf{D}(f,S)$ is obstructed.
\end{lem}

\begin{pf}
As explained in (\cite{WestonUnob}, Section 5.2), in this case $\pi_p$ is special, which translates on the Galois side to the existence of an unramified character $\chi: G_p \to \bar{K}_\lambda^\times$ (where $K$ is the field of Fourier coefficients of $f$ and $K_\lambda$ is its completion at $\lambda$) such that
\[
\rho_{f,\lambda}|_{G_p} \otimes \bar{K}_\lambda  \simeq \left(  \begin{matrix} \epsilon_\ell \chi & \ast \\ 0 & \chi  \end{matrix} \right),
\]
with the upper right corner ramified. Upon reduction this matrix becomes either
\[
A= \left(  \begin{matrix} \bar{\epsilon}_\ell \bar{\chi} & \nu \\ 0 & \bar{\chi}  \end{matrix} \right) \quad \text{or} \quad B=\left(  \begin{matrix}  \bar{\chi} & \nu \\ 0 & \bar{\epsilon}_\ell \bar{\chi}  \end{matrix} \right)
\]
for some $\nu : G_p \to \bar{k}_\lambda$. We note that by (\cite{WestonUnob}, Lemma 5.1), possibility $B$ can only occur if $p^2 \equiv 1$ (mod $\ell$). 

Let $C=\left(  \begin{smallmatrix}  0 & 1 \\ 0 & 0  \end{smallmatrix} \right)$. One computes 
\[
A C A^{-1} = \left(  \begin{matrix}  0 & \bar{\epsilon}_\ell \\ 0 & 0  \end{matrix} \right)
\quad
\text{ and }
\quad
B C B^{-1} = \left(  \begin{matrix}  0 & (\bar{\epsilon}_\ell)^{-1} \\ 0 & 0  \end{matrix} \right),
\]
and so 
\[
(\bar{\epsilon}_\ell \otimes \mathrm{ad} \ \bar{\rho}_{f,\lambda}) \cdot C = \left(  \begin{matrix}  0 & \bar{\epsilon}_\ell^j \\ 0 & 0  \end{matrix} \right)
\]
for $j=0$ or $2$, so $C \in H^0(G_p, \bar{\epsilon}_\ell \otimes \mathrm{ad} \ \bar{\rho}_{f,\lambda})$. If $j=0$, this is obvious; if $j=2$, this follows from the facts that $G_p$ is topologically generated by Frob$_p$, $\epsilon_\ell(\mathrm{Frob}_p)=p$, and $p^2 \equiv 1$ (mod $\ell$). So in either case $H^0(G_p, \bar{\epsilon}_\ell \otimes \mathrm{ad} \ \bar{\rho}_{f,\lambda}) \neq 0$, hence $\mathbf{D}(f,S)$ is obstructed.
\qed \end{pf}

\subsection{Optimal Level Deformation Problems}\label{Thm2Proof}

Let $f \in S_k(\Gamma_1(N),\omega)$ and $g \in S_k(\Gamma_1(N'))$ be newforms in $\mathcal{H}(\bar{\rho})$ with $f$ of optimal level and $N' > N$; by Proposition \ref{Mform}, $N \mid N'$. Let $S$ (resp. $S'$) be the set of places of $\Q$ dividing $N \infty$ (resp. $N' \infty$), so $S \subset S'$.

Write $f = \sum a_n q^n$. Let $K$ be a field containing the Fourier coefficients of both $f$ and $g$, and let $\lambda$ be a prime of $K$ over $\ell$ such that $f \equiv g$ (mod $\lambda$) and hence $\bar{\rho} \simeq \bar{\rho}_{f,\lambda} \simeq \bar{\rho}_{g,\lambda}$. Write $k_\lambda$ for the residue field of $\lambda$. 

Using this notation, we are now ready to prove Theorem \ref{MainThm2}. For the reader's convenience, we restate the theorem.

\noindent \textbf{Theorem 6. }
If $\mathbf{D}(f)$ is unobstructed, then $f$ is of optimal level for $\mathcal{H}(\bar{\rho})$.

\begin{pf} 
We will prove the contrapositive. Keeping the notation from the beginning of Section \ref{Thm2Proof}, let $f$ and $g$ be newforms in $\mathcal{H}(\rhobar)$, with $f$ of optimal level and $g$ of non-optimal level. We will show that $\mathbf{D}(g)$ is obstructed. 

If $\mathbf{D}(f)$ is obstructed, then as noted earlier, this implies $\mathbf{D}(g)$ is also obstructed. So in proving the theorem, we may assume that $\mathbf{D}(f)$ is unobstructed.

We consider separately the primes $p \in S'$ which appear in cases $(1)$, $(2)$, and $(3)$ of Proposition \ref{Mform}. Note that we have an equivalence of deformation problems $\mathbf{D}(g,S')=\mathbf{D}(f,S')$. We write $\mathbf{D}_{\ell}$ for these equivalent deformation problems.

First, suppose $p \mid N'$ is as in case $(3)$, so in particular $p \equiv 1$ (mod $\ell$). Then by Lemma \ref{LemmaPhi} we see that $\mathbf{D}_{\ell}$ is obstructed.

Next, suppose $p \mid N'$ is as in case $(1)$ of the proposition, so $p$ is a prime such that $p \nmid N \ell$, $\alpha(p)=1$, and $\displaystyle pa_p^2 \equiv (1+p)^2 \omega(p)p^{k-1}$ (mod $\lambda$), or equivalently (since $p$ is invertible in $\barFl$),
\[
a_p^2 \equiv (p+1)^2 p^{k-2} \omega(p) \ (\text{mod} \lambda).
\]
Then by Lemma \ref{LemmaAP} we see that $\mathbf{D}_{\ell}$ is obstructed.

Finally, suppose $p \mid N'$ is as in case $(2)$, so $p \equiv -1$ (mod $\ell$) and one of the following holds:
\begin{enumerate}[(a)]
\item $p \nmid N$, $a_p \equiv 0$ (mod $\lambda$), and $\alpha(p)=2$; or
\item $p \mid \mid N$, det $\rho$ is unramified at $p$, and $\alpha(p)=1$.
\end{enumerate}

If $p$ were as in case (b), then actually $p \in S$, and Lemma \ref{LemmaSpecial} shows that $\mathbf{D}(f,S)$ is obstructed. This contradicts our hypothesis on $\mathbf{D}(f)$, so we can ignore this case.

Finally, we must consider case (a), so that $p \equiv -1$ (mod $\ell$), $p \nmid N$, $a_p \equiv 0$ (mod $\lambda$), and $\alpha(p)=2$. Recalling that $\mathbf{D}=\mathbf{D}_{\ell}(f,S')$, Lemma \ref{LemmaAP} gives the obstruction since $a_p \equiv (p+1) \equiv 0$ (mod $\lambda$).
\qed \end{pf}

\begin{rmk}
It is not the case that every minimal, optimal level deformation problem is unobstructed. Indeed, for any prime $p$ we have
\[
\epad \simeq \epb \oplus (\epadzero)
\]
and
\[
H^0(G_p,\epb)\neq 0 \Leftrightarrow p \equiv 1 \quad (\text{mod} \ \ell),
\]
so condition (3) of Theorem \ref{MainThm2} is sharp. Let $\ell=5$, $p=11$, and $k=3$. The space $S_3(\Gamma_1(11),3)$ contains one newform defined over $\Q$ and four newforms which are Galois conjugates defined over $\Q(\alpha)$, where $\alpha$ is a root of $x^4+5x^3+15x^2+15x+5$. The minimal set $S$ for any of these newforms is $S= \{ 11,\infty \} $. Since $S_3(\Gamma_1(1))$ is empty, all of these newforms are of optimal level for their respective mod $\ell$ representations, but since $p \equiv 1$ mod $\ell$ their minimal deformation problems are obstructed.
\end{rmk}

\begin{rmk}
The techniques in this paper cannot rule out the possibility that two (or more) congruent modular forms of optimal level can exist for an unobstructed modular deformation problem.
\end{rmk}

Combining this result with Weston's result (\cite{WestonUnob}, Theorem 1), we have the following corollary.

\begin{cor}
Let $f$ be a newform of level $N$ and weight $k \geq 2$. For infinitely many primes $\ell$, $f$ represents an optimal modular realization of a mod $\ell$ representation $\bar{\rho}: G_{\Q} \to \mathrm{GL}_2 (\bar{\mathbf{F}}_\ell)$. 
\end{cor}

\begin{pf} 
For infinitely many such $\ell$, $\mathbf{D}(f)$ is unobstructed (by Weston), and by Theorem \ref{MainThm2}, this implies that $f$ is of optimal level among modular forms realizing $\rhobar$.
\qed \end{pf}

\begin{rmk}
Actually, there is a much simpler proof of this fact: If $f$ is of nonoptimal level for its mod $\ell$ representation, then there is a modular form $g$ of lower level such that $f \equiv g$. But such a congruence can occur for only finitely many primes $\ell$, which follows from the q-expansion principle and the fact that these spaces of modular forms are finite dimensional.
\end{rmk}

We also get another corollary. For any integer $N$, let $d(N)$ be the number of prime divisors of $N$, i.e. $d(N)=\displaystyle\sum_{p \mid N} 1$.

\begin{cor}
Fix a prime $\lambda$ with residue field $k$ of characteristic $\ell > 3$, suppose $\bar{\rho}: G_\Q \to \mathrm{GL}_2 (k)$ is a modular mod $\lambda$ representation of prime-to-$\ell$ conductor $N$, and let $f$ be a newform of level $N$ such that $\bar{\rho} \simeq \bar{\rho}_{f,\lambda}$. If $g$ is a newform of level $N'$ such that $f \equiv g$ (mod $\lambda$), then
\[
\mathrm{dim}_k H^2(G_{\Q,S},\mathrm{ad} \ \bar{\rho}) \geq d\left( N'/N \right)
\]
where $S$ is a finite set of places containing the primes which divide $N' \infty$.
\end{cor}
\begin{pf}
The proof of Theorem \ref{MainThm2} shows that if $g$ is of nonoptimal level $N'$, then for every prime $p$ dividing $N'/N$, we have $H^0(G_p,\epad) \neq 0$. By equation (\ref{h2dim}) we have
\[
\mathrm{dim}_k H^2(G_{\Q,S},\mathrm{ad}\ \bar{\rho}) \geq \sum_{p \in S} \mathrm{dim}_k H^0(G_p),\epad)
\]
and the corollary follows.

\qed \end{pf}

\section{Explicit Computations in the Supercuspidal Case}




\begin{example}

Let $k=2$, $\ell=11$, and $p=7$. Consider the CM elliptic curve $E$ with Cremona label $49a1$; it is given by 
\[
E: y^2 + xy = x^3 - x^2 -2x-1,
\]
and its associated modular form $f \in S_2(\Gamma_0(49))$ has $q$-expansion
\[
f=q+q^2-q^4-3q^8-3q^9+ \cdots.
\]
The mod $\ell$ Galois representation $\rhobar_{f,\ell}$ is irreducible, and one checks that none of conditions (1)--(8) of Theorem \ref{MainThm1Body} are satisfied. Using Sage \cite{sage}, one also verifies that $\ell \notin$ Cong$(f)$ and so $f$ is of optimal level for this representation. 

Loeffler and Weinstein have incorporated their results from \cite{LW} into the Local Components package of \cite{sage}. Using this, one discovers that $\pi_p$ is supercuspidal. However, since $p^4 \equiv 3$ (mod $\ell$), condition (10) is also not satisfied. Thus $\mathbf{D}(f)$ is unobstructed. 
\end{example}

\begin{example}

This example shows that condition (10) of Theorem \ref{MainThm1Body} is necessary but not sufficient for producing local obstructions at supercuspidal primes. Let us fix $k=3$, $\ell=5$, and $p=7$. Note that $p^2 \equiv -1$ mod $\ell$ and $p^4 \equiv 1$ mod $\ell$. 

Using \cite{sage}, one finds a newform $f$ in $S_3(\Gamma_1(49))$ with a $q$-expansion that begins
\[
f=q + \left(-\frac{1}{92}\alpha^3 + \frac{5}{92}\alpha^2 - \frac{41}{91}\alpha + \frac{229}{92} \right)q^2 +
\left(-\frac{1}{184}\alpha^3 + \frac{5}{184}\alpha^2 - \frac{133}{184}\alpha + \frac{229}{184} \right)q^3 + \cdots
\]

Using the Local Components package of \cite{sage}, one discovers that $\pi_p$ is supercuspidal. Let $E=\Q_p(s)$ be the unramified quadratic extension of $\Q_p$. Let $L=K(\beta)$ where $\beta$ satisfies the polynomial
$\displaystyle
x^2 + \left(  \frac{3}{1288}\alpha^3 + \frac{11}{184}\alpha^2 - \frac{153}{1288} \alpha + \frac{467}{184} \right)x-1.
$ 
Then the character $\chi$ associated to $\pi_p$ is characterized by
\[
\chi: E \to L
\]
\[
s \mapsto \beta, \quad 7 \mapsto 7.
\]
(Here we are viewing $\chi$ as a character of $E^\times$ instead of $G_E$ via local class field theory.) Let $\lambda$ be either of the two primes of $L$ which lies over $\ell$. Then using \cite{sage}, one verifies that $\chi$ and its conjugate $\chi^c$ are equivalent mod $\lambda$ by checking that $\beta - \beta^c$ has positive $\lambda$-valuation. In the notation of Proposition \ref{SCProp}, this shows that $\bar{\psi}=1$; the induction of this character is a symmetric representation, and so $\epb \otimes \mathrm{Ind}_{G_E}^{G_p}\bar{\psi}$ is an invariant $G_p$-representation, hence $H^0(G_p,\epad)=0$. 
\end{example}


\section*{Acknowledgments}

It is a pleasure to thank my thesis advisor, Tom Weston, for suggesting that I study this topic, as well as for teaching me so much mathematics over the years. I am also very grateful to Keenan Kidwell and Ravi Ramakrishna for their many helpful comments on an earlier version of this paper.

\end{document}